\documentclass[11pt]{article}

\def\fpd#1#2{{\displaystyle\frac{\partial #1}{\partial #2}}}
\def\spd#1#2#3{{\displaystyle\frac{\partial^2 #1}
{\partial #2\partial #3}}}
\def\lie#1{{\cal L}_{#1}}
\def\vf#1{\frac{\partial}{\partial #1}}

\def\clift#1{#1^{\scriptscriptstyle{\mathrm{C}}}}
\def\hlift#1{#1^{\scriptscriptstyle{\mathrm{H}}}}
\def\vlift#1{#1^{\scriptscriptstyle{\mathrm{V}}}}

\def\R{\mathcal{R}}

\def\mech#1{#1^{\scriptstyle{\mathrm{m}}}}

\def\onehalf{{\textstyle\frac12}}

\font\frak=eufm10 
\def\goth #1{\hbox{{\frak #1}}}
\def\g{\goth{g}}

\def\cinfty#1{C^{\scriptscriptstyle\infty}(#1)}

\def\pM{{\pi}^{\scriptscriptstyle M}}
\def\pTM{\pi^{\scriptscriptstyle TM}}
\def\psiM{\psi^{\scriptscriptstyle M}}

\def\psiTM{\psi^{\scriptscriptstyle TM}}
\def\Ad{\mathop{\mathrm{ad}}\nolimits}

\def\conn#1#2#3{\setbox1=\hbox{$\scriptstyle{#2}{#3}$}%
\setbox2=\hbox to\wd1{$\hfil\scriptstyle{#1}\hfil$}
\Gamma^{\!\box2}_{\!\box1}}
\def\barconn#1#2#3{\setbox1=\hbox{$\scriptstyle{#2}{#3}$}%
\setbox2=\hbox to\wd1{$\hfil\scriptstyle{#1}\hfil$}
\check{\Gamma}^{\!\box2}_{\!\box1}}

\begin{document}

\title{Invariant Lagrangians, mechanical connections and the
Lagrange-Poincar\'{e} equations}

\author{T.\ Mestdag${}^{a,b}$ and M.\ Crampin${}^{b}$\\[2mm]
{\small ${}^a$ Department of Mathematics, University of Michigan}\\
{\small 530 Church Street, Ann Arbor, MI 48109, USA}\\[1mm]
{\small ${}^b$Department of Mathematical Physics and Astronomy, Ghent University}\\
{\small Krijgslaan 281, B-9000 Ghent, Belgium}}

\date{}

\maketitle

{\small {\bf Abstract.} We deal with Lagrangian systems that are
invariant under the action of a symmetry group.  The mechanical
connection is a principal connection that is associated to Lagrangians
which have a kinetic energy function that is defined by a Riemannian
metric.  In this paper we extend this notion to arbitrary Lagrangians.
We then derive the reduced Lagrange-Poincar\'e equations in a new
fashion and we show how solutions of the Euler-Lagrange equations can
be reconstructed with the help of the mechanical connection.  Illustrative examples confirm the theory.  \\[2mm]
}

\section{Introduction}

The enormous importance of the concept of symmetry in a great number
of applications in physics is beyond any doubt.  Symmetry properties
of mechanical systems in particular have been studied intensively in
the last decades (see e.g.\ the survey in Chapter~3 of the recent
monograph \cite{stagesbook}).  The bulk of the literature, however,
concentrates on the Hamiltonian description of symmetric systems in
which mainly the theory of Poisson manifolds plays an important role.
Less well-known is the process of symmetry reduction for Lagrangian
systems.  When the Lagrangian is invariant under the action of a Lie
group, so that we are dealing with a symmetry group of the
Euler-Lagrange equations, the equations of motion can be reduced to a
new set of equations with fewer unknowns.

In the literature, there are in fact a lot of different paths that
lead to different Lagrangian reduction theories.  For example, if
the configuration space is either a Lie group or a semi-direct
product, the distinct reduction method is called either
Euler-Poincar\'e reduction or semi-direct product reduction. Another
path is the following. In rough terminology, the invariance of the
Lagrangian leads via the Noether theorem to a set of conserved
quantities (the momenta). Whether or not one wants to take these
conserved quantities into account in this reduction process leads to
either the Routh or the Lagrange-Poincar\'e reduction method.  For
more details and some comments on the history of the above mentioned
reduction theories, see e.g.\ \cite{stagesbook} and \cite{MRS}.

In this paper, we concentrate solely on Lagrange-Poincar\'e
reduction, and the associated reduced equations, the so-called
Lagrange-Poincar\'e equations.  The geometric framework for these
equations that has been developed in e.g.\ \cite{Cendra,MRbook}
relies heavily on methods coming from the calculus of variations. It
is shown there that Hamilton's principle can be reduced to a
principle on the reduced space, and that the equations of motion
that follow from this reduced principle are exactly the
Lagrange-Poincar\'e equations. In this paper, we wish to take an
alternative view of Lagrange-Poincar\'e reduction.  We will not take
recourse to variations, but rather interpret the Euler-Lagrange
equations as defining the integral curves of an associated
second-order differential equation field on the velocity phase
space.  We will show that the Lagrange-Poincar\'e equations can be
derived in a relatively straightforward fashion from the
Euler-Lagrange equations, by choosing a suitable adapted frame, or
equivalently by employing well-chosen quasi-velocities.

Next to investigating aspects of reduction, we also focus on the
inverse process, that of reconstruction.  With the aid of a
principal connection, for any invariant system one can reconstruct a
solution of the original problem from a reduced solution.  In the
case of a so-called simple mechanical system (with a Lagrangian of
the form $T-V$ where the kinetic energy $T$ is derived from a
Riemannian metric), it is well-known that such a connection is
naturally available, and it is therefore called `the mechanical
connection'.  We will show in this paper that we can define for any
arbitrary Lagrangian system (under a natural regularity assumption)
a generalized mechanical connection, and that we can employ this
connection in the process of reconstructing integral curves of the
Euler-Lagrange field. The generalized mechanical connection is new.
The construction to be found in Lewis's paper \cite{Lewis}, while
superficially similar is in fact rather different. We emphasize that
the generalized mechanical connection is designed for use in the
reconstruction process only and does not play a role in the
reduction step.

The paper is organised as follows. In the next section we explain
our approach to Lagrangian systems using adapted frames and
quasi-velocities. We extend the notion of a mechanical connection to
arbitrary Lagrangian systems in Section~3. In Section~4 we derive
the Lagrange-Poincar\'e equations within our framework and in
Section~5 we discuss the matter of reconstruction. We end the paper
with some illustrative examples.

\section{Preliminaries}
Consider a Lagrangian function $L$, defined on the tangent bundle
$\tau: TM\to M$ of some manifold $M$ (the configuration space).  In
terms of local coordinates $(x^\alpha,u^\alpha)$ (where the $x^\alpha$
are coordinates on $M$ and the $u^\alpha$ the corresponding fibre
coordinates on $TM$) the Euler-Lagrange equations of $L$ are
\[
\frac{d}{dt}\left(\fpd{L}{{u} ^\alpha}\right)-\fpd{L}{x^\alpha}=0.
\]
We will assume that the Lagrangian is regular, which is to say that
its Hessian with respect to the fibre coordinates,
\[
\spd{L}{{u}^\alpha}{{u}^\beta},
\]
considered as a symmetric matrix, is everywhere non-singular.
In case the Lagrangian is regular, the Euler-Lagrange equations can
be written explicitly in the form of a system of differential
equations $\ddot{x}^\alpha=\Gamma^\alpha(x,\dot{x})$. These
equations may be thought of as defining a vector field $\Gamma$ on
$TM$, a second-order differential equation field, namely
\[
\Gamma={u}^\alpha\vf{x^\alpha}+\Gamma^\alpha\vf{{u}^\alpha};
\]
we call this the Euler-Lagrange field of $L$.  In fact when $L$ is
regular the Euler-Lagrange field is uniquely determined by the fact
that it is a second-order differential equation field and satisfies
\[
\Gamma\left(\fpd{L}{{u}^\alpha}\right)-\fpd{L}{x^\alpha}=0.
\]

In this paper we will work with regular Lagrangians and their
Euler-Lagrange fields, determined by the version of the Euler-Lagrange
equations given immediately above.  However, we will often not use
coordinates; instead, we will express everything in terms of certain
convenient local bases of vector fields.  We must therefore explain
how to modify the formulation above to take account of this
difference.

There are two canonical ways of lifting a vector field $Z$ on $M$ to
one on $TM$.  Firstly, the vector field on $TM$ with the property that
its flow is tangent to the flow of $Z$ on $M$ is called the complete
lift of $Z$, and is denoted by $\clift{Z}$.  In terms of coordinates
\[
\clift Z=Z^\alpha \vf{x^\alpha}+\fpd{Z^\beta}{x^\alpha}
{u}^\alpha\vf{{u}^\beta},\quad \mbox{where }Z=Z^\alpha(x)\fpd{}{x^\alpha}.
\]
The second canonical way of lifting the vector field $Z$
from $M$ to $TM$ is the vertical lift:
\[
\vlift{Z}= Z^\alpha\fpd{}{{u}^\alpha};
\]
$\vlift{Z}$ is tangent to the fibres of the projection $\tau: TM \to
M$, and on the fibre over $m\in M$ it coincides with the constant
vector field $Z_m$.  In general, we have $T\tau\circ\clift{Z}=Z\circ\tau$
while $T\tau\circ\vlift{Z}=0$ (regarding the vector fields as sections
of the appropriate bundles).  The Lie brackets of complete and vertical
lifts of vector fields $Y$ and $Z$ on $M$ are given by $[\clift Y,\clift Z] =
[Y,Z\clift]$, $[\clift Y,\vlift Z] = [Y,Z\vlift]$ and $[\vlift
Y,\vlift Z]=0$. It is worth emphasising that although the map
$Z\mapsto\vlift{Z}$ is $C^\infty(M)$-linear, the map
$Z\mapsto\clift{Z}$ is only $\R$-linear:\ in fact for $f\in
C^\infty(M)$
\[
\clift{(fZ)}=f\clift{Z}+\dot{f}\vlift{Z}
\]
where $\dot{f}$ is the so-called total derivative of $f$, given in
coordinates by
\[
\dot{f}=u^\alpha\fpd{f}{x^\alpha}
\]
(a function on $TM$).  More details on all of this can be found in
e.g.\ \cite{CP}.

Now take a local basis of vector fields on $M$, say $\{Z_\alpha\}$.
Then $\{\clift{Z_\alpha},\vlift{Z_\alpha}\}$ is a local basis of vector
fields on $TM$.  Moreover, we can introduce new fibre coordinates
$v^\alpha$ on $TM$, adapted to the local basis $\{Z_\alpha\}$, as
follows: for any tangent vector $v_m\in T_mM$, the $v^\alpha$ are the
components of $v$ with respect to the basis $Z_\alpha(m)$, that is,
$v_m=v^\alpha Z_\alpha(m)$. Such fibre coordinates are called quasi-velocities.

From the coordinate expressions of the vertical and complete lifts
above it is easy to see that if $\{Z_\alpha\}$ is such a basis then
the equations
\begin{equation}
\Gamma(\vlift{Z_\alpha}(L))-\clift{Z_\alpha}(L)=0
\end{equation}
are equivalent to the Euler-Lagrange equations. Moreover, the fact that
$\Gamma$ is a second-order differential equation field means that it
takes the form
\begin{equation}
\Gamma= v^\alpha\clift{Z_\alpha}+\hat{\Gamma}^\alpha\vlift{Z_\alpha}
\end{equation}
where the $v^\alpha$ are the quasi-velocities corresponding
to the basis $\{Z_\alpha\}$.

We will also need a coordinate-independent version of the Hessian.
The Hessian of $L$ at $u\in TM$ is the symmetric bilinear form $g$ on
$T_mM$, $m=\tau(u)$, given by $g(v,w)=\vlift{v}(\vlift{w}(L))$, where
the vertical lifts are to $u$. Alternatively we can regard $g$ as
operating on pairs of vector fields $Y$, $Z$ on $M$:\ then $g(Y,Z)$ is a
function on $TM$, depending bilinearly on its
arguments, such that $g(Y,Z)=\vlift{Y}(\vlift{Z}(L))$; the fact that it
is symmetric in $Y$ and $Z$ is a consequence of the bracket relation
$[\vlift{Y},\vlift{Z}]=0$ stated earlier.

We turn now to the question of symmetries and invariance.

Throughout the paper we will assume that the configuration space $M$
of the Lagrangian system comes equipped with a free and proper left
action $\psiM: G\times M\to M$ of a group $G$, which eventually will
be the symmetry group of the Lagrangian system under consideration.
In taking the group to act to the left we follow the convention of
\cite{MRbook}; changing to a right action (as is the standard in e.g.\
\cite{KN}) has the effect only of a sign change in some of our results
and formulae.  The projection $\pM: M\to M/G$ on the set of
equivalence classes gives $M$ the structure of a principal fibre
bundle.  The action on $M$ induces an action $\psiTM_g = T\psiM_g$ on
the tangent manifold $TM$, so $\pTM:TM\to TM/G$ is also a principal
fibre bundle.

A tensor field $t$ on $M$ is invariant under the action of $G$ if for
all $g\in G$,
\[
t_{\psiM_g m}(\psi^{T^*M}_g\alpha_m,\ldots, \psiTM_g v_m,\ldots)
=t_m(\alpha_m,\ldots,v_m,\ldots), \qquad\alpha_m\in T^*_mM,  v_m\in T_mM,
\]
where $\langle \psi^{T^*M}_g\alpha_m, \psiTM_g v_{m} \rangle =\langle
\alpha_m, v_m\rangle$.  For convenience we will assume that the group
$G$ is connected.  Then the invariance property is equivalent to the
infinitesimal condition $\lie{\tilde\xi}t=0$ for all $\xi\in\g$ (the
Lie algebra of $G$), where $\tilde\xi$ on $M$ is the fundamental
vector field associated with $\xi$.  The proof of the equivalence
relies on the fact that the flow of $\tilde\xi$ on $M$ is $\psi_t: m
\mapsto \psiM_{\exp(t\xi)}m$.  In particular, a function $f$ on $M$ is
invariant if and only if it satisfies $\tilde\xi(f)=0$, and a vector
field $X$ on $M$ is invariant if and only if $[{\tilde \xi},X]=0$, for
all $\xi\in\g$.

A tensor field on $TM$ is invariant if and only if all of its Lie
derivatives by fundamental vector fields of the induced action on $TM$
vanish.  The flow of such a vector field is given by $v_m\mapsto
T\psiM_{\exp(t\xi)}v_m$ for some $\xi\in\g$, and is therefore tangent
to $m\to \psiM_{\exp(t\xi)}m $, the flow of the vector field
$\tilde\xi$ on $M$.  That is to say, the fundamental vector field
corresponding to $\xi\in\g$ of the action of $G$ on $TM$ is just
$\clift{\tilde\xi}$.

For the remainder of the paper we will suppose that $L$ is a
Lagrangian function on $TM$ invariant under the induced action of a
connected group $G$ on $TM$, so that $L(\psiTM_gv)=L(v)$, or
equivalently $\clift{\tilde \xi}(L)=0$ for all $\xi\in\g$.

We will work with a local basis $\{{\tilde E}_a,X_i\}$ of vector
fields on $M$ which is adapted to the bundle structure $M\to M/G$ in
the following way.  The vector fields ${\tilde E}_a$ are the
fundamental vector fields corresponding to a basis $\{E_a\}$ of $\g$;
the $X_i$ are the $G$-invariant horizontal lifts of a coordinate basis
of vector fields on $M/G$, where the horizontal lift is associated
with some principal connection $\omega$ on $M\to M/G$.  The
quasi-velocities corresponding to such a basis will be denoted by
$(v^a,v^i)$, so that for every tangent vector $v_m$, $v_m=v^a{\tilde
E}_a(m)+v^iX_i(m)$.

It will also be convenient to have a basis $\{{\hat E}_a,X_i\}$ that
consists only of invariant vector fields.  Let $U\subset M/G$ be an
open set over which $M$ is locally trivial.  The projection $\pM$ is
locally given by projection onto the first factor in $U\times G \to
U$, and the (left) action of $G$ by $\psiM_g(x,h)=(x,gh)$.  The vector
fields on $M$ defined by
\[
{\hat E}_a: (x,g) \mapsto \widetilde{(\Ad_{g} E_a)}(x,g) =
\psiTM_g\big(\tilde{E}_a(x,e)\big)
\]
are invariant.  The relation between the sets $\{{\hat E}_a\}$ and
$\{{\tilde E}_a\}$ can be expressed as ${\hat E}_a(x,g)=
A_a^b(g){\tilde E}_b(x,g)$ where $(A_a^b(g))$ is the matrix
representing $\Ad_{g}$ with respect to the basis $\{E_a\}$ of $\g$.
In particular, $A^b_a(e)=\delta^b_a$.  The quasi-velocities
corresponding to the basis $\{{\hat E}_,X_i\}$ will be denoted by
$(w^a,v^i)$, where in fact $v^a=A^a_bw^b$.

The left-invariant vector fields $\hat{E}_a$ satisfy the bracket
relations $[{\hat E}_a,{\hat E}_b] = C_{ab}^c {\hat E}_c$ where the
$C^c_{ab}$ are the structure constants of the Lie algebra (so that
the Lie algebra bracket satisfies $[E_a,E_b]=C_{ab}^c E_c$).  For
the fundamental vector fields, on the other hand, we have $[{\tilde
E}_a,{\tilde E}_b] = -C_{ab}^c {\tilde E}_c$.  Other Lie brackets of
basis vector fields are:\ $[X_i,{\tilde E}_a]=0$ by invariance,
whence $[X_i,{\hat E}_a]= X_i(A_a^c)\bar{A}_c^b{\hat
E}_b=\Upsilon_{ia}^b {\hat E}_b$ say, where $(\bar{A}_a^b)$ is the
matrix inverse to $A^b_a$; and $[X_i,X_j] = K^a_{ij}{\hat E}_a$,
where the $K^a_{ij}$ are the components of the curvature of
$\omega$, regarded as a $\g$-valued tensor field.

By expressing the fact that the vector fields ${\hat E}_a$ are
invariant in the form $[{\tilde E}_b, {\hat E}_a] = 0$ we find that
\[
{\tilde E}_b(A^c_a)-A^d_aC^c_{bd}=0.
\]

With respect to the basis $\{{\hat E}_a,X_i\}$ the Euler-Lagrange
field is of the form
\begin{equation}
\Gamma = w^a \clift{\hat E}_a + v^i \clift{ X}_i + \Gamma^a
\vlift{\hat E}_a + \Gamma^i \vlift{X}_i.
\end{equation}
The functions $\Gamma^i$ and $\Gamma^a$ can be determined from the
Euler-Lagrange equations
\begin{eqnarray}
\Gamma(\vlift{X_i}(L))-\clift{X_i}(L)&=&0\nonumber\\
\Gamma(\vlift{\hat{E}_b}(L))-\clift{\hat{E}_b}(L)&=&0.
\end{eqnarray}

We will now deduce from the assumed invariance of $L$ that the
Euler-Lagrange field $\Gamma$ is invariant, using this machinery.
That is to say, we will show that if $\clift{\tilde E}_a(L)=0$ for
all $a$, then $[\clift{\tilde E}_a,\Gamma]=0$ also.  One easily
verifies that the quasi-velocities corresponding to an invariant
basis are invariant (and in any case the derivation is given in full
later); it follows that
\[
[\clift{\tilde E}_b,\Gamma] = \clift{\tilde E}_b(\Gamma^i)
\vlift{X_i} + \clift{\tilde E}_b(\Gamma^a) \vlift{{\hat E}_a}.
\]
From the above Euler-Lagrange equations it easily follows that
\begin{eqnarray*}
0&=& \clift{\tilde E}_b(\Gamma(\vlift{X_i}(L)))-
\clift{\tilde E}_b (\clift{X_i}(L))\\
&=& [\clift{\tilde E}_b,\Gamma](\vlift{X_i}(L)) +
\Gamma ( \clift{\tilde E}_b(\vlift{X_i}(L)) )
- [\clift{\tilde E}_b, \clift{X_i}] (L) -
\clift{X_i}(\clift{\tilde E}_b (L)) \\
&=&[\clift{\tilde E}_b,\Gamma](\vlift{X_i}(L))
+ \Gamma ( [\clift{\tilde E}_b,
\vlift{X_i}](L)) + \Gamma (\vlift{X_i} (\clift{\tilde E}_b(L)))\\
&=&
[\clift{\tilde E}_b , \Gamma](\vlift{X_i}(L)).
\end{eqnarray*}
Likewise, $[\clift{\tilde E}_b,\Gamma](\vlift{{\hat E}_c}(L))=0$. From
the expression for $[\clift{\tilde E}_b,\Gamma]$ we obtain
\[
\clift{\tilde E}_b(\Gamma^i)g(X_i,X_j)
+ \clift{\tilde E}_b(\Gamma^a)g(\hat{E}_a,X_j)=0=
\clift{\tilde E}_b(\Gamma^i)g(X_i,\hat{E}_c) +
\clift{\tilde E}_b(\Gamma^a)g(\hat{E}_a,\hat{E}_c).
\]
Now at any point $u$ of $TM$ we may regard $[\clift{\tilde
E}_b,\Gamma]$ as the vertical lift of some vector $w\in T_mM$,
$m=\tau(u)$.  From the last displayed equations we conclude that
$g(w,X_j(m))=0=g(w,\hat{E}_c(m))$.  The required conclusion that $w=0$
follows from the assumed non-singularity of $g$.  Notice that it
follows that $\clift{\tilde E}_b(\Gamma^i)=0$ and $\clift{\tilde
E}_b(\Gamma^a)= 0$.

\section{The mechanical connection}

A connection on $\pM:M\to M/G$ is a left splitting $\omega$ of the
short exact sequence
\[
0\to M\times \g \to TM \stackrel{T\pi^M}{\to} (\pM)^* T(M/G)\to 0.
\]
We identify $M\times \g$ as a subbundle of $TM$ by means of
$(m,\xi)\mapsto \tilde\xi(m)$.  We will use the symbol $\omega$ for
the following two related objects: $\omega$ may be thought of as a
type $(1,1)$ tensor field on $M$, for which, in particular,
$\omega(\tilde\xi)=\tilde\xi$.  On the other hand, we will also use
$\omega$ for its projection onto $\g$; we then have
$\omega(\tilde\xi(m)) = \xi$.  With the second interpretation, if
$\omega$ satisfies $\omega(\psiTM_g v)= \Ad_{g}\omega(v)$ the
connection is said to be principal.  The infinitesimal version of this
invariance property is that $\lie{\tilde\xi}\omega =0$, for all $\xi$,
where $\omega$ is now interpreted as a $(1,1)$-tensor on $M$.  The
kernel of $\omega$ (in either interpretation) defines a distribution
on $TM$, invariant when the connection is principal, which is called
the horizontal distribution.  We can also define the connection by
specifying its horizontal distribution.

In the case of a so-called simple mechanical system, one can
associate a principal connection with the system in a natural way. In
this case the Lagrangian takes the form $L=T-V$ where $T$ is a
kinetic energy function, defined by a Riemannian metric $g$ on $M$,
and $V$ is a function on $M$, the potential energy. The symmetry
group of the Lagrangian consists of those isometries of $g$ which
leave $V$ invariant. The mechanical connection $\omega$ can be
defined by taking for its horizontal subspaces the orthogonal
complements of the tangent spaces to the fibres of $\pM$.

The construction of the mechanical connection relies heavily on the
availability of a Riemannian metric.  Now the components of the metric
with respect to a coordinate basis coincide with those of the Hessian
of the simple Lagrangian with respect to the fibre coordinates.
However, the Hessian of an arbitrary Lagrangian does not provide a
Riemannian metric on configuration space.  Nevertheless, the notion of
a mechanical connection may be generalized in such a way that the
Hessian plays the same role in relation to the generalized mechanical
connection as the metric does to the mechanical connection for a
simple system.  This we now show.

Recall that the Hessian of $L$ at $u\in TM$ is the symmetric bilinear
form $g$ on $T_mM$, $m=\tau(u)$, given by
$g(v,w)=\vlift{v}\vlift{w}(L)$.  As such, it can be interpreted as a
tensor field along the tangent bundle projection $\tau: TM \to M$, as
we show in the following paragraphs.

A vector field along $\tau$ is a section of the pullback
bundle $\tau^*TM \to TM$. Such a section can in an equivalent way be
interpreted as a map $X: TM\to TM$ with the property that
$\tau\circ X=\tau$. A vector field along $\tau$ takes the
local form
\[
X = X^\alpha(x,u) \fpd{}{x^\alpha}.
\]
A vector field $Z$ on $M$ can be interpreted as the vector field
$Z\circ\tau$ along $\tau$.  We will call such vector fields along
$\tau$ `basic', and use the same symbol for the vector field on $M$
and the basic vector field along $\tau$.

By taking for 1-forms along $\tau$ the sections of $\tau^*T^*M \to
TM$, we can obtain in the usual manner a $\cinfty{TM}$-module of
tensor fields along $\tau$. (For more information see e.g.\
\cite{MCS,Szilasi}.) In particular, we can interpret the Hessian of
a Lagrangian as the symmetric (0,2) tensor field along $\tau$ given
by
\[
g= g_{\alpha\beta }(x,u) dx^\alpha \otimes dx^\beta,\qquad
g_{\alpha\beta} = \frac{\partial^2 L}{\partial {u}^\alpha\partial
{u}^\beta}.
\]
Then as we pointed out earlier, if $X,Y$ are vector fields on $M$ then
$g(X,Y)=\vlift X (\vlift Y (L))$.

We show now that when $L$ is invariant, $g$ is also invariant in an
appropriate sense.

The action of the symmetry group $G$ on $TM$ induces an action of $G$ on
$\tau^*TM$ by $(g,(u_m,v_m)) \mapsto (\psiTM_g u_m,\psiTM_g v_m)$ for
$g\in G$.  Then, for example, a vector field $X$ along
$\tau$ is invariant if
\[
X(\psiTM_g v) = \psiTM_g (X(v)), \qquad \forall g\in G.
\]
The infinitesimal version of this property will follow from a more
general construction in \cite{MikeDavid}, which we adapt
to the current context.

Let $\varphi:TM\to TM$ be fibred over $f:M\to M$. Then we can extend
$\varphi$ to a map $\bar{\varphi}:\tau^*(TM)\to \tau^*(TM)$, fibred
over $\varphi$, as follows:\ for $u_m,v_m\in TM$, set
\[
\bar{\varphi}(u_m,v_m)=(\varphi(u_m),Tf(v_m)),
\]
where $Tf:TM\to TM$ is the tangent map to $f$.  We have $\tau\circ
\varphi=\tau\circ Tf=f$, so $\tilde{\varphi}$ is well defined.  In
coordinates $(x^\alpha,u^\alpha,v^\alpha)$, with
$\varphi(x,u)=(f^\alpha(x),\varphi^\alpha(x,u))$, we have
\[
\bar{\varphi}(x,u,v)
=\left(f^\alpha(x),\varphi^\alpha(x,u),\fpd{f^\alpha}{x^\beta}v^\beta\right).
\]

Given any vector field $X$ along $\tau$ and any fibred
diffeomorphism $\varphi:TM\to TM$ we can define a new vector field
along $\tau$, $\varphi_\sharp X$, by $\varphi_\sharp
X(v)=\bar{\varphi}(X(\varphi^{-1}(v)))$. We can thus define a Lie
derivative operator of vector fields on $M$ on the set of vector
fields along $\tau$:\ for any vector field $Z$ on $M$ and vector
field $X$ along $\tau$ set
\[
\lie{Z} X=\frac{d}{dt}(\psi_{(-t)\sharp} X)_{t=0},
\]
where $\psi_t$ is the flow of $\clift Z$. If we think of $X$ as a
derivation from functions on $M$ to functions on $TM$, we have for a
function $F$ on $M$
\[
\lie{Z} X(F)=\clift Z( X(F))-X(Z(F)).
\]
In coordinates
\[
\lie{Z}\left(X^\alpha\vf{x^\alpha}\right)=
\left(Z^\beta\fpd{X^\alpha}{x^\beta}+
\fpd{Z^\beta}{x^\gamma}{u}^\gamma\fpd{X^\alpha}{u^\beta}
-X^\beta\fpd{Z^\alpha}{x^\beta}\right)\vf{x^\alpha}.
\]
If $X$ is a basic vector field along $\tau$ then $\lie{Z} X$ is the
basic vector field $[Z,X]$ along $\tau$, where the bracket is the
bracket of vector fields on $M$.

The operator  $\lie{Z}$ has all the usual properties of the Lie
derivative, except that for functions $F$ on $TM$,
$\lie{Z}(FX)=\clift{Z}(F)X+F\lie{Z}X$.  So if we define the action
of $\lie{Z}$ on functions by $\lie{Z}F=\clift{Z}(F)$, we can extend
the action to the tensor algebra of $\tau^*TM$ in the usual way. In
particular, for a type $(0,2)$ tensor field $g$ and vector fields
$X$, $Y$ along $\tau$,
\[
(\lie{Z}g)(X,Y)=\clift{Z}(g(X,Y))-g(\lie{Z}X,Y)-g(X,\lie{Z}Y).
\]

Consider the case where $X$ is an invariant vector field along
$\tau$, and where $Z$ is a fundamental vector field $\tilde\xi$.
The flow of $\tilde\xi$ is exactly $f_t:m \mapsto
\psiM_{\exp(t\xi)}m$. By definition of a complete lift, the flow of
$\clift{\tilde\xi}$ is the tangent of the flow of $\tilde\xi$, that
is  $\psi_t: v \mapsto T\psiM_{\exp(t\xi)}v =\psiTM_{\exp(t\xi)}v$.
Since $(\psi_{-t})^{-1}=\psi_t$, we find that
\[ ((\psi_{-t})_\sharp X) (v) = {\bar\psi}_{-t} \Big(X(\psi_t v\Big) =
\Big(\psi_{-t}(\psi_t v), Tf_{-t}(X(\psi_t v ))\Big) = (v,X(v)),
\]
where we have used the invariance property of $X$ with respect to
group elements of the form $\exp(t\xi)$.  So $(\psi_{-t})_\sharp
X=X$ and we conclude that $\lie{\tilde\xi} X=0$ for all $\xi\in\g$.
Since $G$ is supposed to be connected, a standard argument shows that
this criterion is in fact sufficient.

The vertical lift operation described earlier can easily be extended to
vector fields along $\tau$. From the coordinate expression of the
Lie derivative it is easy to see that $(\lie{Z}X\vlift)= [\clift
Z,\vlift X]$. It follows in particular that a vector field $X$ along
$\tau$ is invariant if and only if its vertical lift is invariant as
a vector field on $TM$.

With a similar argument to the one for vector fields along $\tau$,
we can conclude that a tensor field $g$ along $\tau$ is invariant if
and only if $\lie{\tilde\xi}g=0$ for all $\xi\in\g$. In fact it is
enough to show that this tensor vanishes when its arguments are
basic vector fields $X,Y$ along $\tau$. Let $g$ be the Hessian of an
invariant Lagrangian. Recall that for $X$, $Y$ basic,
$g(X,Y)=\vlift{X}(\vlift{Y}(L))$. We have
\begin{eqnarray*}
(\lie{\tilde{\xi}}g)(X,Y)&=& \clift{\tilde{\xi}}(g( X, Y))
-g(\lie{\tilde{\xi}} X, Y)
-g( X,\lie{\tilde{\xi}} Y)\\
&=&\clift{\tilde{\xi}}(\vlift{X}\vlift{Y}(L))
-\vlift{[\tilde{\xi},X]}\vlift{Y}(L)
-\vlift{X}\vlift{[\tilde{\xi},Y]}(L)\\
&=&[\clift{\tilde{\xi}},\vlift{X}]\vlift{Y}(L)
+\vlift{X}\clift{\tilde{\xi}}\vlift{Y}(L)
-\vlift{[\tilde{\xi},X]}\vlift{Y}(L)
-\vlift{X}\vlift{[\tilde{\xi},Y]}(L)\\
&=&[\clift{\tilde{\xi}},\vlift{X}]\vlift{Y}(L)
-\vlift{[\tilde{\xi},X]}\vlift{Y}(L)
+\vlift{X}[\clift{\tilde{\xi}},\vlift{Y}](L)
-\vlift{X}\vlift{[\tilde{\xi},Y]}(L)\\
&=&0,
\end{eqnarray*}
since by assumption $\clift{\tilde{\xi}}(L)=0$.

We turn now to the definition of the generalized mechanical connection.  We
first make some remarks about connections in this general
context.

The short exact sequence $0\to M\times\g\to TM\to {\pM}^*T(M/G)\to 0$
extends in a natural way to an exact sequence of vector bundles over
$TM$,
\[
0\to TM\times\g \to \tau^*TM \to ({\pM}\circ\tau)^* T(M/G)\to 0,
\]
where the first space is spanned by basic vector fields ${\tilde\xi}$
along $\tau$.  We will call an invariant left splitting of this
sequence an invariant connection along $\tau$.  Equivalently, an
invariant connection along $\tau$ is a (1,1)-tensor field $\omega$
along $\tau$ with the property that $\lie{\tilde\xi}\omega=0$ for all
$\xi\in\g$.  For a basic vector field $\tilde\eta$ along $\tau$ we
have $\omega(\tilde\eta)=\tilde\eta$.  A vector field $X$ along $\tau$
which satisfies $\omega(X)=0$ is said to be horizontal.

From the fact that $\omega(\tilde\eta)=\tilde\eta$ it follows
automatically that
\[
(\lie{\tilde\xi}\omega)(\tilde\eta)= \omega
([\tilde\xi,\tilde\eta]) - [\tilde\xi,\omega(\tilde\eta)]= 0.
\]
So, for invariance, the only condition we need to check is that for all
horizontal vector fields $X$ along $\tau$,  $\lie{\tilde\xi}X$
is horizontal also. In fact
\[
(\lie{\tilde\xi}\omega)(X)=
\omega(\lie{\tilde\xi}X) -\lie{\tilde\xi}(\omega(X))=
\omega(\lie{\tilde\xi}X),
\]
so $\lie{\tilde\xi}\omega=0$ if and only $\lie{\tilde\xi}X$ is
horizontal whenever $X$ is horizontal.

Any principal connection on $M\to M/G$ can be extended to
an invariant connection along $\tau$ in a trivial way.

In the case of an invariant Lagrangian system we can define a
connection along $\tau$ by adapting the construction of the mechanical
connection of a simple mechanical system; the result is the
generalized mechanical connection, which we denote by $\mech\omega$.
A vector field $X$ along $\tau$ is horizontal with respect to
$\mech\omega$ if it is in the orthogonal complement of $TM\times\g$
with respect to $g$, that is if
\begin{equation}
g({\tilde \xi},X)=0, \qquad \forall \xi\in\g.
\end{equation}
Of course, this defines a splitting only if we suppose that $g$ is
non-singular when restricted to the set of fundamental vector
fields; this will be the case in particular if $g$ is everywhere
positive-definite.

We now show that $\mech\omega$ is invariant, by showing that if $X$
is a horizontal vector field along $\tau$, so is $\lie{\tilde\xi}X$
for all $\xi\in\g$.  For a horizontal $X$
\begin{eqnarray*}0&=& (\lie{\tilde\xi}g)(\tilde\eta,X) =
\clift{\tilde\xi}(g(\tilde\eta,X)) -
g([\tilde\xi,\tilde\eta],X) - g(\tilde\eta, \lie{\tilde\xi}X)\\
&=&   g (\widetilde{[ \xi, \eta]},X)- g(\tilde\eta,\lie{\tilde\xi}X)\\
&=& - g(\tilde\eta,\lie{\tilde\xi}X),
\end{eqnarray*}
so $\lie{\tilde\xi}X$ is indeed horizontal, as claimed.

We can lift $\mech\omega$ in a natural way to define a principal
connection $\mech\Omega$ on $TM\to TM/G$. The process by which we
obtain $\mech\Omega$ from $\mech\omega$ is a variant of the so-called
vertical lift of a principal connection $\omega$ on $\pM:M\to M/G$ to $TM\to
TM/G$: the vertical lift of $\omega$, considered here as a
$\g$-valued 1-form on $M$, is just $\tau^*\omega$, a $\g$-valued
1-form on $TM$. It is easy to show that $\tau^*\omega$ is a principal
connection on $\pTM:TM\to TM/G$:\ see for example \cite{Paper1}.

The vertical lift of the generalized mechanical connection
$\mech\omega$ is the $\g$-valued 1-form $\mech\Omega$ on $\pTM$ defined by
\begin{equation}
\mech\Omega(W) = \mech\omega(\tau_*W)
\end{equation}
for any vector field $W$ on $TM$; here $\tau_*W$ is the projection of
$W$ regarded as a vector field along $\tau$.  Thus a vector field $W$
on $TM$ is horizontal if $g(\tau_*W,\tilde{\eta})=0$ for all
$\eta\in\g$, where as usual $\tilde{\eta}$ is regarded as a basic
vector field along $\tau$.  It is clear that the horizontal subspace
of $T_uTM$ is complementary to the tangent to the fibre of $\pTM$,
provided that $g_u$ is non-singular on the tangent to the fibre of
$\pM$, as before.  Note in particular the paradoxical-sounding fact
that every vertical vector field is horizontal:\ to be precise, every
vector field on $TM$ which is vertical with respect to the tangent
bundle projection $\tau$ is horizontal with respect to the connection
$\mech\Omega$.  Now $\mech\Omega$ defines a principal connection if
$[\clift{\tilde\xi},W]$ is horizontal whenever $W$ is horizontal.
From the definition of the generalized Lie derivative above it
is easy to see that for any vector field $W$ on $TM$ we have
\[
\lie{Z}(\tau_*W)=\tau_*[\clift Z,W].
\]
Then, indeed,
\[
g(\tau_*[\clift{\tilde{\xi}},W],\tilde{\eta})=
g(\lie{{\tilde{\xi}}}(\tau_*W),\tilde{\eta})
=-g(\tau_*W,\lie{{\tilde{\xi}}}\tilde{\eta})=
-g(\tau_*W,[\tilde{\xi},\tilde{\eta}])=0.
\]
This confirms that $\mech\Omega$ is a principal connection on $\pTM$.

The generalized mechanical connection $\mech\omega$ and its vertical
lift $\mech\Omega$ are intrinsic to the invariant Lagrangian system;
that is to say, it is not necessary to invoke a principal connection
on $\pM: M\to M/G$ to define them. Nevertheless it is often convenient to
work with a local basis of vector fields $\{\tilde{E}_a,X_i\}$ on
$M$, with $X_i$ invariant, as described earlier; and this implicitly
involves a principal connection on $\omega$ on $\pM$. In general
$\omega$ will not be directly related to $\mech\omega$; but for a
simple mechanical system we can take for $\omega$ the mechanical
connection, in which case $\mech\omega$ is the natural extension of
$\omega$ to an invariant connection along $\tau$, and $\mech\Omega$
is the vertical lift of $\omega$.

We will end this section by expressing the generalized mechanical
connection and its vertical lift in terms of the basis
$\{\tilde{E}_a,X_i\}$.

Let us express the components of the Hessian $g$ in terms of the
basis $\{{\tilde E}_a,X_i\}$ as follows:
\[
g(\tilde{E}_a,\tilde{E}_b)=g_{ab},\quad g(X_i,X_j)=g_{ij},\quad
g(X_i,\tilde{E}_a)=g_{ia}=g_{ai}=g(\tilde{E}_a,X_i)
\]
(in general these will be functions on $TM$, not $M$, of course).  We
have $g_{ab}=\vlift{\tilde E}_a(\vlift{\tilde E}_b (L))$, etc.  Recall
that to define the generalized mechanical connection we have assumed
that $(g_{ab})$ is non-singular. If we set
\[
{\bar X}_i=X_i-g^{ab}g_{ib}{\tilde E}_a =X_i+B^a_i {\tilde E}_a
\]
then the $\bar{X}_i$ are vector fields along $\tau$ which are
horizontal with respect to the generalized mechanical connection, and
so
\[
\mech\omega({\tilde E}_a) = {\tilde E}_a, \qquad
\mech\omega({\bar X}_i)=0.
\]
The invariance of the Hessian $g$ amounts for its coefficients
to
\[
\clift{\tilde{E}_a}(g_{ij})=0,\quad
\clift{\tilde{E}_a}(g_{bc})=C^d_{ab}g_{cd}+C^d_{ac}g_{bd},\quad
\clift{\tilde{E}_a}(g_{ic})=C^d_{ac}g_{id}.
\]
It follows that $\clift{\tilde E}_a(B^b_i) = B^c_i C^b_{ca}$.  It is
now easy to see that $\bar{X}_i$ is an invariant vector field along
$\tau$:
\[
\lie{{\tilde E}_a} {\bar X}_i =
\lie{{\tilde E}_a} {X}_i +\clift{\tilde E}_a (B^b_i){\tilde E}_b +
B^b_i \lie{{\tilde E}_a}{\tilde E}_b =
[{\tilde E}_a, {X}_i] + C^b_{ca} B^c_i {\tilde E}_b +
B^b_i [{\tilde E}_a,{\tilde E}_b]=0.
\]
Furthermore, let us define vector fields
\[
\clift{\bar X}_i =\clift{X_i} + B^a_i \clift{{\tilde E}_a}
\]
on $TM$.  (The notation $\clift{\bar X}_i$ is not intended to imply
that these vector fields are complete lifts of vector fields on $M$.)
One can easily verify that the vector fields $\clift{\bar X}_i$ are
invariant:\ $[\clift{\tilde E}_a, \clift{\bar X}_i]=0$.  For the
lifted connection $\mech\Omega$ we get
\[
\mech\Omega(\clift{{\tilde E}_a})= \clift{{\tilde E}_a},
\qquad\mech\Omega(\vlift{{\tilde E}_a}) = 0, \qquad
\mech\Omega({\clift{\bar X}}_i)=0 \quad \mbox{and}\quad
\mech\Omega(\vlift{X}_i) = 0.
\]

\section{The reduced Euler-Lagrange field}

As before, we assume that the Lagrangian $L$ is invariant under a
symmetry group $G$, so that $\clift{\tilde \xi}(L)=0$.  Then $L$
defines a function $l$ on $TM/G$, the reduced Lagrangian, such that
$L=l\circ\pTM$, where $\pTM$ is the projection of the principal fibre
bundle $TM\to TM/G$.  We showed above that the Euler-Lagrange field
$\Gamma$ is also invariant.  As an invariant vector field, $\Gamma$ on
$TM$ defines a $\pTM$-related reduced vector field $\check \Gamma$ on
$TM/G$:\ due to the invariance of $\Gamma$, the relation $T\pTM
\big(\Gamma (v)\big)= \check\Gamma \big(\pTM(v)\big)$ is independent
of the choice of $v\in TM$ within the equivalence class of $\pTM(v)\in
TM/G$.

Our aim in this section is to give an explicit expression for the
reduced Euler-Lagrange field $\check\Gamma$ on $TM/G$ in terms of the
reduced Lagrangian $l$.

From now on we will use local coordinates on $M$ defined as follows.
Let $U\subset M/G$ be an open set over which $M$ is locally trivial,
so that $(\pM)^{-1}(U)\simeq U\times G$, and which is a coordinate
neighbourhood in $M/G$.  We take coordinates $(x^\alpha)=(x^i,x^a)$ on
a suitable open subset of $(\pM)^{-1}(U)$ (containing $U\times{e}$)
such that $(x^i)$ are coordinates on $U$, $(x^a)$ coordinates on the
fibre $G$.  Then $\tilde{E}_a$ and $\hat{E}_a$ can be identified with
vector fields on $G$, right- and left-invariant respectively, and so
are independent of the $x^i$.  Recall that $X_i$ is an invariant
vector field on $M$ projecting onto $\partial/\partial x^i$ on $M/G$.
If we set
\[
{X}_i = \fpd{}{x^i} - \gamma_i^b(x^i,x^a) {\hat E}_b
\]
then invariance of the $X_i$ means that
\[
\fpd{\gamma^b_i}{x^a}=0.
\]
Recall that we have set $[X_i,\hat{E}_a]=\Upsilon^b_{ia}\hat{E}_b$.
Then $\Upsilon_{ia}^b =-\gamma_i^cC^b_{ca}$.

An invariant vector field $W$ on $TM$ projects onto a vector field
$\check W$ on $TM/G$, so that (as sections)
\[
T\pTM\circ W=\check{W}\circ\pTM.
\]
Furthermore, if $F$ is an invariant function on $TM$ and if $f$ is its
reduction to a function on $TM/G$ then
\[
W(F) = W(f\circ\pTM)=\check W(f) \circ\pTM.
\]
We will now show in detail that the coordinate functions $x^i$, $v^i$
and $w^a$ (where $w^a$ and $v^i$ are the quasi-velocities
corresponding to our chosen invariant basis $\{{\hat E}_a,X_i\}$) are
invariant functions on $TM$.  As a consequence, they can be used as
coordinates on $TM/G$.  Then the action of $W$ on these functions
completely determines $\check W$.  We are in particular interested in
the reduced fields of the invariant vector fields $\clift{X}_i$,
$\clift{\bar X}_i$, $\vlift{X}_i$ and $\vlift{\hat E}_a$ on $TM$.

The following observation can easily be verified. For any vector
field $Z$, function $f$ and 1-form $\theta$ on $M$,
\[
\clift Z (f) = Z(f), \quad \vlift Z(f)=0,\quad
\clift{Z}(\vec{\theta})=\overrightarrow{\lie{Z}\theta}, \quad
\vlift{Z}(\vec{\theta})=\tau^*\theta(Z),
\]
where $\vec\theta$ stands for the fibre-linear function on $TM$
defined by the 1-form $\theta$. Note that if $\{Z_\alpha\}$ is a
local basis of vector fields on $M$ and $\theta^\alpha$ the dual
local basis of 1-forms, then the quasi-velocities $v^\alpha$
corresponding to the given vector-field basis, regarded  as functions
on $TM$, are given by $v^\alpha=\vec{\theta}^\alpha$.

Let $\{\varpi^a,\vartheta^i\}$ be the dual 1-form basis of the basis
$\{{\hat E}_a,X_i\}$. Then $\vec{\varpi}^a=w^a$ and
$\vec{\vartheta}^i=v^i$. We have, for example,
\begin{eqnarray*}
(\lie{X_i}\varpi^a)(\hat{E}_b)&=&
X_i(\delta^a_b)-\varpi^a([X_i,\hat{E}_b])= -\Upsilon_{ib}^a  \\
(\lie{X_i}\varpi^a)(X_j)&=&-\varpi^a([X_i,X_j])=-K^a_{ij},
\end{eqnarray*}
so $\clift{X_i}(w^a)=-K^a_{ij}v^j - \Upsilon_{ib}^aw^b$. The
relevant derivatives, obtained in the case of the quasi-velocities
by similar calculations, are
\[
\begin{array}{lllll}
\clift{X_i}(x^j)=\delta_i^j, &\quad&
\clift{X_i}(v^j)=0,&\quad&
\clift{X_i}(w^a)=-K^a_{ij}v^j -\Upsilon_{ib}^aw^b,\\
\vlift{X_i}(x^j)=0, &&
\vlift{X_i}(v^j)=\delta^j_i,&&
\vlift{X_i}(w^a)=0,\\
\clift{{\hat E}_a}(x^j)=0, &&
\clift{\hat{E}_a}(v^i)=0,&&
\clift{\hat{E}_a}(w^b)= \Upsilon^b_{ia}v^i+C_{ac}^bw^c, \\
\vlift{{\hat E}_a}(x^j)=0,&&
\vlift{\hat{E}_a}(v^i)=0,&&
\vlift{\hat{E}_a}(w^b)=\delta^b_a,\\
\clift{{\tilde E}_a}(x^j)=0, &&
\clift{\tilde{E}_a}(v^i)=0,&&
\clift{\tilde{E}_a}(w^b)=0,\\
\vlift{{\tilde E}_a}(x^j)=0,&&
\vlift{\tilde{E}_a}(v^i)=0, &&
\vlift{\tilde{E}_a}(w^b)={\bar{ A}}_a^b.
\end{array}
\]
From the penultimate row we conclude that $x^i$, $v^i$ and $w^a$ are
invariant. Therefore
\begin{eqnarray*} &&
T\pTM\circ \clift{ X_i}  = \left(\fpd{}{x^i} + (-K^a_{ij}v^j -
\Upsilon_{ib}^aw^b)\fpd{}{w^a}\right) \circ\pTM =
T\pTM \circ\clift{{\bar X}_i} ,\\
&& T\pTM\circ \vlift{X_i}  = \fpd{}{v^j}\circ\pTM, \\
&& T\pTM\circ\clift{\hat E}_a=
(\Upsilon^b_{ia}v^i+C_{ac}^bw^c)\fpd{}{w^b}\circ\pTM,\qquad\qquad
T\pTM\circ \vlift{\hat E}_a =\fpd{}{w^a}\circ\pTM.
\end{eqnarray*}

As we noted earlier the Euler-Lagrange field is invariant.  The
Euler-Lagrange equations, in terms of the invariant basis
$\{\hat{E}_a,X_i\}$, are
\begin{eqnarray*}
\Gamma(\vlift{X_i}(L))-\clift{X_i}(L)&=&0\\
\Gamma(\vlift{\hat{E}_b}(L))-\clift{\hat{E}_b}(L)&=&0.
\end{eqnarray*}
Using invariance we see that
\[
\Gamma(\vlift{X_i}(L))=\Gamma(\vlift{\check{X}_i}(l)\circ\pTM)
=\check{\Gamma}(\vlift{\check{X}_i}(l))\circ\pTM,
\]
and similarly for the other terms in the Euler-Lagrange equations.
These may therefore be written, entirely in terms of vector fields and
functions on $TM/G$, as
\begin{eqnarray*}
\check{\Gamma}(\vlift{\check{X}_i}(l))-\clift{\check{X}_i}(l)&=&0\\
\check{\Gamma}(\vlift{\check{E}_b}(l))-\clift{\check{E}_b}(l)&=&0,
\end{eqnarray*}
where $\vlift{\check{E}_b}$ and $\clift{\check{E}_b}$ are the
projections onto $TM/G$ of the invariant vector fields
$\vlift{\hat{E}_b}$ and $\clift{\hat{E}_b}$.  We can now employ the
coordinate expressions for the reduced vector fields
$\vlift{\check{X}_i}$ etc.\ obtained earlier to express these
equations as
\begin{eqnarray}
&&\check{\Gamma}\left( \fpd{l}{v^i}\right) -\fpd{l}{x^i} =
(K^a_{ik}v^k + \Upsilon^a_{ib}w^b)\fpd{l}{w^a}\nonumber\\
&&\check{\Gamma}\left( \fpd{l}{w^a}\right) =
(\Upsilon^b_{ia}v^i + C^b_{ac}w^c)\fpd{l}{w^b}.
\end{eqnarray}
Recall that, with respect to the basis $\{X_i,{\hat E}_a\}$, the
Euler-Lagrange field is of the form $ \Gamma = w^a \clift{\hat E}_a
+ v^i \clift{ X}_i + \Gamma^a \vlift{\hat E}_a + \Gamma^i
\vlift{X}_i$. Given that $C_{ab}^c w^aw^b=0$ and $K^a_{ij}v^iv^j=0$,
the coordinate expression of the reduced field is simply
$\check\Gamma=v^i\partial/\partial x^i + \Gamma^i \partial/\partial
v^i + \Gamma^a \partial/\partial w^a$; as we observed earlier,
$\Gamma^i$ and $\Gamma^a$ are invariant, so may be considered as
functions on $TM/G$.  The above equations, usually written in the
form
\begin{eqnarray*}
&&\frac{d}{dt}\left( \fpd{l}{v^i}\right) -\fpd{l}{x^i} =
(K^a_{ik}v^k + \Upsilon^a_{ib}w^b)\fpd{l}{w^a}\\
&&\frac{d}{dt}\left( \fpd{l}{w^a}\right) = (\Upsilon^b_{ia}v^i +
C^b_{ac}w^c)\fpd{l}{w^b},
\end{eqnarray*}
are the so-called Lagrange-Poincar\'e equations (see e.g.\
\cite{Cendra}).

The specific structure of the manifold and the symmetry group at
hand can lead to some interesting subcases. For example, if the
symmetry group $G$ happens to be Abelian, all terms containing the
structure coefficients of the Lie algebra vanish. From the last
equation we then obtain that the `momentum' $\partial l/\partial
w^a$ is constant, let's say $\mu_a$. If moreover the matrix
$(\partial^2 l/\partial w^a\partial w^b)$ is non-singular, the
relation $\partial l/\partial w^a=\mu_a$ can be rewritten in the
form $w^a=\rho^a(x,v)$. With that, the first equation becomes a
second order differential equation in the variables $x^i$. By
introducing Routh's (reduced) function
\[
{\mathcal R}^\mu (x,v)= l(x,v, \rho(x,v)) - \mu_a \rho^a(x,v)
\]
the first equation can equivalently be rewritten as
\[
\frac{d}{dt}\left( \fpd{{\mathcal R}^\mu}{v^i}\right)
-\fpd{{\mathcal R}^\mu}{x^i} = K^a_{ik}v^k\mu_a.
\]
This equation is known as Routh's (reduced) equation for an Abelian
symmetry group. Routh's reduction process can be extended to
non-Abelian symmetry groups. However, the extension requires that
the two steps in the process above (first reduction and then
restriction to a level set of momentum) be alternated. For more
details, see \cite{Routhpaper} and \cite{MRbook}. There is no
obvious way to relate the Lagrange-Poincar\'e equations to the
reduced Routh equations in the case of a non-Abelian symmetry group.

Another interesting case occurs when the manifold $M$ is a product
$N\times G$. Then, we can choose the connection to be trivial and
all terms involving connection coefficients will vanish. Finally, if
the manifold is in fact the Lie group $G$, the equations
\[\frac{d}{dt}\left( \fpd{l}{w^a}\right) = C^b_{ac}w^c\fpd{l}{w^b}
\] are known as the Euler-Poincar\'e equation, see e.g.\
\cite{MRbook}.

\section{Reconstruction}

Now that we have derived the reduced form of the Euler-Lagrange
equations it remains to consider the problem of reconstruction:\
suppose we can find a solution of these equations, that is, an
integral curve of $\check\Gamma$, how do we reconstruct from it a
solution of the original equations, that is, an integral curve of
$\Gamma$?

There is in fact a standard method for reconstructing integral curves
of an invariant vector field from reduced data, which makes use of
connection theory.  Let $\pi:M\to B$ be a principal fibre bundle with
group $G$.  An invariant vector field $X$ on $M$ defines a
$\pi$-related reduced vector field $\check{X}$ on $B$.  Let us
suppose that $M$ is equipped with a principal connection $\omega$.
Let ${\check v}(t)$ be an integral curve of $\check{X}$ (in
$B$).  Let $m$ be a point of $M$ in the fibre over ${\check v}(0)$:\
we aim to find the integral curve of $X$ though $m$.  We first
lift ${\check v}(t)$ to $M$ by using the connection to form its
horizontal lift through $m$, $\hlift{\check{v}}(t)$:\ this is the
unique curve in $M$ projecting onto $\check v$ such that
$\omega(\dot{\hlift{\check v}})=0$ and $\hlift{\check{v}}(0)=m$. Now
let $v(t)$ be the integral curve of $X$ through $m$.
Since $v(t)$ also projects onto $\check{v}(t)$ there is a curve $t\mapsto
g(t)\in G$, with $g(0)=e$, such that
$v(t)=\psiM_{g(t)}\hlift{\check{v}}(t)$.  On differentiating this equation
we obtain
\[
\dot{v} =
\psiTM_{g}\left(\widetilde{\vartheta(\dot{g})}\circ\hlift{\check{v}}
+\dot{\hlift{\check{v}}}\right)
\]
where $\vartheta$ is the Maurer-Cartan form of $G$ (so that
$\vartheta(\dot{g}(t))$ is a curve in $\g$).  But
\[
\dot{v}(t)=X(v(t))=
X(\psiM_{g(t)}\hlift{\check{v}}(t))=\psiTM_{g(t)}X(\hlift{\check{v}}(t)),
\]
since $X$ is invariant, from which it follows that
\[
\widetilde{\vartheta(\dot{g})}\circ\hlift{\check{v}}
+\dot{\hlift{\check{v}}}=X\circ\hlift{\check{v}}.
\]
The first term on the left-hand side is vertical, the second
horizontal, so this equation is simply the decomposition of
$X\circ\hlift{\check{v}}$ into its horizontal and vertical components
with respect to $\omega$.  Thus $g(t)$ must satisfy the so-called
reconstruction equation
\[
\vartheta(\dot{g}) = \omega(X\circ\hlift{\check{v}}),
\]
where $\omega$ is taken to be the connection 1-form, so that the
right-hand side is a curve in $\g$.  This is a differential equation
for the curve $g(t)$, and has a unique solution with specified initial
value (see for example \cite{Sharpe}).  Thus the curve $g(t)$ is
uniquely determined by the equation and the initial condition
$g(0)=e$.  Conversely, if $g(t)$ is the solution of the reconstruction
equation such that $g(0)=e$ then
$v(t)=\psiM_{g(t)}\hlift{\check{v}}(t)$ is the integral curve of $X$
through $m$.

We can use this method to obtain an integral curve of $\Gamma$ from
one of $\check{\Gamma}$ by using the generalized mechanical connection
on $\pTM:TM\to TM/G$. The reconstruction equation in this case is
\begin{equation}
\vartheta(\dot g) = \mech\Omega(\Gamma\circ\hlift{\check v}).
\end{equation}
We have therefore to find the vertical component of $\Gamma$ with
respect to $\mech\Omega$. This is not completely straightforward
because we have expressed $\Gamma$ in terms of the invariant basis to
obtain the reduced Euler-Lagrange equations in the previous section,
whereas $\mech\Omega$, as a $\g$-valued 1-form, is specified by
\[
\mech\Omega(\clift{{\tilde E}_a})= E_a,
\qquad\mech\Omega(\vlift{{\tilde E}_a}) = 0, \qquad
\mech\Omega({\clift{\bar X}}_i)=0, \qquad
\mech\Omega(\vlift{X}_i) = 0,
\]
where
\[
\clift{\bar X}_i =\clift{X_i} + B^a_i \clift{{\tilde E}_a}
=\clift{X_i} -g^{ab}g_{bi} \clift{{\tilde E}_a}.
\]
Recall that $\hat{E}_a=A_a^b\tilde{E}_b$. From the general
properties $(fZ\clift) = f \clift Z + \dot f \vlift Z$ and
$(fZ\vlift)=f\vlift Z$ we find that
\[
\clift{\hat E}_a = A_a^b \clift{\tilde E}_b + \dot{A}_a^b
\vlift{\tilde E}_b , \qquad \vlift{\hat E}_a = A_a^b \vlift{\tilde
E}_b,
\]
and therefore
\[
\mech\Omega(\clift{\hat{E}_a})=A_a^bE_b,\qquad
\mech\Omega(\vlift{\hat{E}_a})=0.
\]
Moreover, $\mech\Omega(\clift{ X}_i)=B^a_iE_a$. Thus with
\[
\Gamma = w^a \clift{\hat E}_a + v^i \clift{ X}_i + \Gamma^a
\vlift{\hat E}_a + \Gamma^i \vlift{X}_i.
\]
we have
\[
\mech\Omega(\Gamma)=(A^a_bw^b-B^a_iv^i)E_a.
\]
The reconstruction equation is therefore
\[
\vartheta(\dot{g}(t))=
(A^a_b(\hlift{\check{v}}(t))w^b(t)-B^a_i(\hlift{\check{v}}(t))v^i(t))E_a;
\]
since $w^a$ and $v^i$ are invariant, $w^b(t)$ and $v^i(t)$ are just
their values on $\check{v}(t)$.  In fact in the coordinate system we
used in the previous section, corresponding to a chart $U\times G$
on $M$, the (left) action on $M$ is simply given by $\psiM_g(x,h) =
(x,gh)$, and the induced action on $TM$ by $\psiM_g (x,h,w^a,v^i) =
(x,gh,w^a,v^i)$ (to use a somewhat bastardized but self-explanatory
notation).  We assume that we are able to calculate the integral
curve $\check v(t)=(x(t),w^a(t),v^i(t))$ of the reduced vector field
$\check\Gamma$ through $\pTM(v_0)$ for some point $v_0\in TM$.  The
horizontal lift of $\check v$ is a curve in $TM$ of the form
$\hlift{\check{v}}(t)=(x(t),h(t),w^a(t),v^i(t))$, where $h(t)$ is a
curve in $G$ to be determined by $v_0$ and by the relation
$\mech\Omega({\dot{\hlift{\check v}}})=0$.  Moreover, since we have
identified the fibres of $\pM:M\to M/G$ with $G$, $A^a_b$ is
effectively a function on $G$.  The reconstruction equation can
therefore be written
\begin{equation}
\vartheta(\dot{g}(t))=
(A^a_b(h(t))w^b(t)-B^a_i(\hlift{\check{v}}(t))v^i(t))E_a,
\end{equation} and the integral curve of $\Gamma$ is just
$v(t)=(x(t),g(t)h(t),w^a(t),v^i(t))$.

It may be of interest to express $\Gamma$ in its vertical and
horizontal components with respect to $\mech\Omega$. In the first
place,
\[
\Gamma=(A^a_bw^b-B^a_iv^i)\clift{\tilde{E}_a}+v^i\clift{\bar{X}_i}
+(\Gamma^a+w^b\dot{A}^c_b\bar{A}^a_c)\vlift{\hat{E}_a}+\Gamma^i\vlift{X}_i,
\]
where $(\bar{A}^b_a)$ is the matrix inverse to $(A^b_a)$.  To proceed
further we need a more revealing expression for the term involving
$\dot{A}^c_b$, the total derivative of $A^c_b$.  We can rewrite
$\dot{A}^c_b$ as
\[
w^d\hat{E}_d(A^c_b)+v^iX_i(A^c_b).
\]
Now
\[
\hat{E}_d(A^c_b)=A_d^e\tilde{E}_e(A^c_b)=A_d^eC_{ef}^cA^f_b,
\]
and therefore
\[
w^bw^d\hat{E}_d(A^c_b)=C_{ef}^c(A_d^ew^d)(A^f_bw^b)=0.
\]
Recall that $X_i(A^c_b)=\Upsilon_{ib}^dA^c_d$. Thus
\[
w^b\dot{A}^c_b\bar{A}^a_c=\Upsilon_{ib}^av^iw^b,
\]
and the coefficient of $\vlift{\hat{E}_a}$ is
$\Gamma^a+\Upsilon_{ib}^av^iw^b$.

By taking into account the fact that $K^a_{ij}v^iv^j=0$, together
with the expression just obtained, it easily follows that
the horizontal part of $\Gamma$ projects onto $\check\Gamma$, as it
must.

The above connection is not the same as the one we have used in
\cite{Paper1} for the reconstruction of an arbitrary (not
necessarily Lagrangian) second-order field $\Gamma$ on $TM$.  There
we started with an arbitrary principal connection $\omega$ on $M$
and formed its vertical lift $\Omega$ on $TM$, as we described
earlier.  For projectable vector fields $W$ on $TM$ we have
$\Omega(W)=\omega(\tau_*W)$ (regarding connections as type $(1,1)$
tensor fields).  Thus
\[
\Omega(\clift{{\tilde E}_a})= \clift{{\tilde E}_a},\qquad
\Omega(\vlift{{\tilde E}_a}) = 0, \qquad
\Omega(\clift{X}_i)=0 \quad \mbox{and}\quad
\Omega(\vlift{X}_i) = 0.
\]
With respect to the vertical lift connection $\Omega$
the reconstruction equation is just
\begin{equation}
\vartheta({\dot g}) =
A^b_a(h(t))w^b(t) E_b,
\end{equation}
where the curve $h(t)$ in $G$ is now of course determined by the
horizontal lift of $\check{v}$ with respect to $\Omega$, not
$\mech\Omega$ as before.

Note that in the case of a simple mechanical system, where we can take
$\omega$ to be the mechanical connection, the vertical lifts $\Omega$
and $\mech\Omega$ coincide.  In that case, the Hessian coincides with
the Riemannian metric, and the horizontal vector fields $X_i$ are
orthogonal to the vector fields ${\tilde E}_a$, from which $g_{ai}=0$
and therefore also $B^a_i=0$.

\section{Illustrative examples}

\subsection{A charged particle in a magnetic field and Wong's equations}

We apply the above introduced machinery to the Kaluza-Klein
formulation of a charged particle in a magnetic field, see e.g.\
\cite{MRbook}. We will consider two steps of abstraction. In the
first step, we assume given a Riemannian manifold on which a group
$G$ acts freely and properly to the left as isometries and we make
the further stipulation that the vertical part of the metric (that
is, its restriction to the fibres of $\pi^M:M\to M/G$) comes from a
bi-invariant metric on $G$. We write down the geodesic equations of
the metric, by interpreting them as the Euler-Lagrange equations for
the kinetic energy Lagrangian.  The reduced equations in such a case
are known as Wong's equations~\cite{Cendra,Mont}. In the second
step, we take the manifold to be $\mathrm{E}^3\times\mathrm{S}$ and
the metric to be
 of Kaluza-Klein type.

We will denote the metric by $g$.  The fact that the symmetry group
acts as isometries means that the fundamental vector fields
$\tilde{\xi}$ are Killing fields:\ $\lie{\tilde{\xi}}g=0$.  It
follows that the components of $g$ with respect to the members of an
invariant basis $\{{\hat E}_a,X_i\}$ are themselves invariant. We
will set $g({\hat E}_a,{\hat E}_b)=h_{ab}$ and $g(X_i,X_j)=g_{ij}$.
We will use the mechanical connection, which means that
$g(\hat{E}_a,X_i)=0$.  Since both $h_{ab}$ and $g_{ij}$ are
$G$-invariant functions, they pass to the quotient; in particular,
the $g_{ij}$ are the components with respect to the coordinate
fields of a metric on $M/G$, the reduced metric.

The further assumption about the vertical part of the metric has the
following implications. It means in the first place that $\lie{{\hat
E}_c}g({\hat E}_a,{\hat E}_b)=0$ (as well as
$\lie{\tilde{E}_c}g({\hat E}_a,{\hat E}_b)=0$), and secondly that
the $h_{ab}$ must be independent of the coordinates $x^i$ on $M/G$,
which is to say that they must be constants.  From the first
condition, taking into account the bracket relations $[{\hat
E}_a,{\hat E}_b]=C^c_{ab}{\hat E}_c$, we easily find that the
$h_{ab}$ must satisfy $h_{ad}C^d_{bc}+h_{bd}C^d_{ac}=0$. Recall that
if we set
\[
X_i=\vf{x^i}-\gamma_i^a{\hat E}_a
\]
for some $G$-invariant coefficients $\gamma_i^a$, we get
$\Upsilon^b_{ia}=\gamma_i^cC_{ac}^b$, and therefore
$h_{ac}\Upsilon^c_{ib}+h_{bc}\Upsilon^c_{ib}=0$.

The geodesic equations may be derived from the Lagrangian
\[
L=\onehalf g_{\alpha\beta}u^\alpha u^\beta =\onehalf
g_{ij}v^iv^j+\onehalf h_{ab}w^aw^b.
\]
It is of course $G$-invariant. We may therefore apply
Lagrange-Poincar\'{e} reduction, which gives the reduced equations
\begin{eqnarray*}
\frac{d}{dt}(g_{ij}v^j)-\onehalf\fpd{g_{jk}}{x^i}v^jv^k&=&
-(K^a_{ij}v^j+\Upsilon^a_{ib}w^b)h_{ac}w^c\\
\frac{d}{dt}(h_{ab}w^b)&=&(\Upsilon^b_{ia}v^i+C^b_{ac}w^c)h_{bd}w^d.
\end{eqnarray*}
Now $\Upsilon^a_{ib}h_{ac}$ is skew-symmetric in $b$ and $c$, and
$C_{ac}^bh_{bd}$ is skew-symmetric in $c$ and $d$, so the final
terms in each equation vanish identically. Let $\conn ijk$ be the
connection coefficients of the Levi-Civita connection of the reduced
metric $g_{ij}$. Then, we may write the equations in the form
\begin{eqnarray*}
g_{ij}\left(\ddot{x}^j+\conn jkl\dot{x}^k\dot{x}^l\right)
&=&-h_{bc}K^c_{ij}\dot{x}^jw^b\\
h_{ab}\left(\dot{w}^b+\Upsilon^b_{ic}\dot{x}^iw^c\right)&=&0,
\end{eqnarray*}
using the skew-symmetry of $\Upsilon^c_{ib}h_{ac}$ again in the
second equation. Given that $K^c_{ij}$ is of course skew-symmetric
in its lower indices, these equations are equivalent to
\begin{eqnarray*}
\ddot{x}^i+\conn ijk\dot{x}^j\dot{x}^k &=&g^{ik}h_{bc}K^c_{jk}\dot{x}^jw^b\\
\dot{w^a}+\Upsilon^a_{jb}\dot{x}^jw^b&=&0.
\end{eqnarray*}
These are Wong's equations.

Let us now take $M$ to be $\mathrm{E}^3\times\mathrm{S}$, with
coordinates $(x^i,\theta)$. Let $A_i$ be the components of a
covector field on $\mathrm{E}^3$, and define a metric $g$ on $M$,
the Kaluza-Klein metric, by
\[
g=\delta_{ij}dx^i\odot dx^j+(A_idx^i+d\theta)^2
\]
where $(\delta_{ij})$ is the Euclidean metric. The Kaluza-Klein
metric admits the Killing field $E=\partial/\partial\theta$. The
vector fields $X_i=\partial/\partial x^i-A_i\partial/\partial\theta$
are orthogonal to $E$ and invariant; moreover
$g_{ij}=g(X_i,X_j)=\delta_{ij}$, while $g(E,E)=1$. Finally
\[
[X_i,X_j]=\left(\fpd{A_i}{x^j}-\fpd{A_j}{x^i}\right)\vf{\theta}.
\]
Putting these values into the reduced equations above we obtain
\[
\ddot{x}^i=w\dot{x}^j\left(\fpd{A_i}{x^j}-\fpd{A_j}{x^i}\right),
\qquad \dot{w}=0.
\]
These are the equations of motion of a particle of unit mass and
charge $w$ in a magnetic field whose vector potential is $A_i dx^i$.

\subsection{A worked-out example}

We will consider the Lie group $G$ of the affine line.  An element of
this group is an affine map ${\bf R}\to{\bf R}: t\mapsto \exp(\theta)
t+ \phi$ and can be represented by the matrix
\[
\left(\begin{array}{ll} \exp\theta & \phi \\ 0 & 1 \end{array}\right).
\]
The identity element is just $t\mapsto t$ (the identity matrix) and
multiplication on the left of $(\theta_2,\phi_2)$ by
$(\theta_1,\phi_1)$ is given by the composition of the two affine
maps, i.e.\ the element
\[
(\theta_1,\phi_1) *(\theta_2,\phi_2)= (\theta_1+\theta_2,
\exp(\theta_1)\phi_2+\phi_1).
\]
 The corresponding Lie algebra is given by
the set of matrices of the form
\[
\left( \begin{array}{ll} a & b \\ 0 & 0 \end{array}\right).
\]

The manifold $M$ of interest is $G\times{\bf R}$.  The action on the
manifold is given by left translation on the $G$ factor of $G\times{\bf
R}$.  We will write $x^0=x$ for the coordinate on ${\bf R}$.  For this
action and this manifold there is a trivial principal connection
$\omega$, with $\gamma^0_a=0$.

A basis of fundamental vector fields is
\[
{\tilde E}_1 = \fpd{}{\theta} + \phi\fpd{}{\phi}, \qquad {\tilde
E}_2 = \fpd{}{\phi}.
\]
The vector fields that are horizontal with respect to the trivial connection
all lie in the direction of $X=\partial/\partial x$. The adapted
coordinates $(v^i,v^a)$ are therefore $v^0={\dot x}$ and $v^1 = \dot\theta$,
$v^2= \dot\phi-\phi\dot\theta$.

The Lie algebra bracket is given by
$[{\tilde E}_1,{\tilde E}_2]=-{\tilde E}_2$.  The complete and
vertical lifts of this basis are
\[
\clift{\tilde E}_1 = \fpd{}{\theta} + \phi \fpd{}{\phi} +
\dot\phi\fpd{}{\dot\phi}, \qquad \clift{\tilde E}_2 =  \fpd{}{\phi}
 ,\qquad \vlift{\tilde E}_1 = \fpd{}{{\dot \theta}} + \phi\fpd{}{\dot\phi}, \qquad
\vlift{\tilde E}_2 = \fpd{}{{\dot \phi}}.
\]

An invariant basis of vector fields is given by
$\{{\hat E}_1,{\hat E}_2,X\}$, where
\[
{\hat E}_1 = \fpd{}{\theta}, \qquad {\hat E}_2 =
\exp(\theta)\fpd{}{\phi},
\]
and the coordinates with respect to this basis are $v^0=\dot x$,
$w_1={\dot \theta}$ and $w_2= \exp(-\theta){\dot \phi}$. The
complete and vertical lifts of the above basis are
\[
\clift X = \fpd{}{x}, \quad \vlift X = \fpd{}{\dot x},\quad
\clift{\hat E}_1 = \fpd{}{\theta}, \quad \clift{\hat E}_2 =
\exp(\theta)\Big(\fpd{}{\phi} +{\dot \theta}\fpd{}{{\dot
\phi}}\Big),\quad \vlift{\hat E}_1 = \fpd{}{{\dot \theta}}, \quad
\vlift{\hat E}_2 = \exp(\theta)\fpd{}{{\dot \phi}}.
\]

Finally, the matrix ${\mathcal A}$, defined by the relation ${\hat
E}_a(x,g)= {A}_a^b(g){\tilde E}_b(x,g)$, is here
\[
{A}(g) = \left( \begin{array}{cc} 1 & 0 \\ -\phi & \exp(\theta)
\end{array}\right).
\]
At the identity of the Lie group, the matrix ${A}$ is the identity
matrix, as it should be.

If we use the invariant fibre coordinates $(v^0,w^a)$, the induced
action on $TM$ is simply\\
$\psi^{TM}_{(\phi_1,\theta_1)}(x,(\phi,\theta), \dot x, w_1, w_2) =
(x,(\phi_1,\theta_1)*(\phi,\theta), \dot x, w_1, w_2)$. Since the
coordinates $(x, \dot x, w_1, w_2)$ can be interpreted as
coordinates on $TM/G= T{\bf R} \times TG/G = T{\bf R}\times \g$,
invariance of the Lagrangian simply means that the group variables
do not explicitly appear in the Lagrangian, when it is written in
terms of the invariant fibre coordinates.

We will work with the Lagrangian
\[
L = \onehalf{\dot\theta}^2  + q\dot x\dot\theta + \onehalf {\dot x}^2
+ \ln(\exp(-\theta)\dot\phi) ,
\]
where $q$ is a constant. The Lagrangian is clearly not of the simple
type. One easily verifies that the Lagrangian is invariant:
$\clift{\tilde E}_1 (L)=0= \clift{\tilde E}_2(L)$. In the invariant
fibre coordinates $(w_i=(w_1,w_2),\dot x)$, the Lagrangian is
\[
L=\onehalf w_1^2  + q\dot x w_1+\onehalf{\dot x}^2 + \ln(w_2),
\]
so, indeed, $(\phi,\theta)$ do not appear explicitly.  This Lagrangian
is also invariant under the obvious ${\bf R}$-action on the manifold,
but we will not take this into consideration.

The Hessian matrix in the basis $\{ {\tilde E}_a, X\}$ is here
\[
g =\left(
\begin{array}{ccc} 1- \frac{\phi^2}{{\dot\phi}^2} &
-\frac{\phi}{{\dot\phi}^2} & q\\
-\frac{\phi}{{\dot\phi}^2} & -\frac{1}{{\dot\phi}^2} & 0 \\
q & 0 & 1
\end{array}\right).
\]
The determinant of $g$ is $(q^2-1 )/{\dot\phi}^2$, so the
Lagrangian is regular as long as $q^2\neq 1$. The upper left (2,2)
matrix represents $(g_{ab})$. It is non-singular since its
determinant is $-1/{\dot\phi}^2$. Its inverse is
\[
(g^{ab}) = \left( \begin{array}{cc} 1&  -\phi \\ -\phi &
\phi^2-{\dot\phi}^2 \end{array} \right).
\]

The vector field ${\bar X}$ along $\tau$ that projects onto
$\partial/\partial x$ on $M/G$ and is horizontal for the generalized
mechanical connection $\mech\omega$ is
\begin{eqnarray*}
{\bar X} &=& \fpd{}{x} - g^{bc}g_{cx}{\tilde E}_b \\&=&\fpd{}{x} - q {\tilde E}_1 +q\phi{\tilde E}_2\\
&=& \fpd{}{x} - q \fpd{}{\theta}.
\end{eqnarray*}
This vector field is in fact a basic vector field along $\tau$.
Although the Lagrangian and the Hessian are not of the simple type,
in this example the generalized mechanical connection turns out to
be derived from a principal connection on $M\to M/G$, and is not of the
most general case of an invariant connection on the pullback bundle.

The corresponding vector field that is horizontal with respect to the
connection $\mech\Omega$ is
\begin{eqnarray*}
\clift{\bar X} &=& \fpd{}{x} - q \clift{\tilde E}_1 +
q\phi\clift{\tilde E}_2
\\ &=& \fpd{}{x} - q \fpd{}{\theta} -q \dot\phi\fpd{}{\dot\phi}.
\end{eqnarray*}

Let us look now at the dynamics. First, let us pretend that we do not
know that the system exhibits symmetry and solve directly the
Euler-Lagrangian equations for this Lagrangian. The equations for
$x$, $\theta$ and $\phi$ are here, respectively:
\[
q\ddot\theta+\ddot x=0,\qquad \ddot\theta+q\ddot x+1=0,\qquad
-\frac{\ddot\phi}{{\dot\phi}^2}=0.
\]
From the $x$-equation, it is again clear that there is also ${\bf
R}$-symmetry. The solution of the system is easy to find. With the
obvious notations for the integration constants, we obtain
\[
x(t) = -\onehalf\frac{qt^2}{q^2-1}+{\dot x}_0t+x_0,\quad
\theta(t)
= \onehalf \frac{t^2}{q^2-1}+{\dot\theta}_0t+\theta_0,\quad
\phi(t)= {\dot\phi}_0t+\phi_0.
\]
We will assume that ${\dot\phi}_0>0$.

We now apply the technique of symmetry reduction and reconstruction.
We first need a solution of the Lagrange-Poincar\'e equations. Since
for the current example the connection coefficients of the trivial
connection vanish, these equations become
\[
\frac{d}{dt}\left(\fpd{l}{w^b}\right) =  \fpd{l}{w^a}C^a_{bd}w^d,
\qquad \frac{d}{dt}\left(\fpd{l}{{\dot x}}\right) - \fpd{l}{x}=0.
\]
The reduced Lagrangian on $TM/G$ is
\[
l(x,\dot x,w_1,w_2)=
\onehalf w_1^2  + q\dot x w_1+\onehalf{\dot x}^2 + \ln(w_2),
\]
and the reduced equations are
\[
\dot{w}_1  +q\ddot x  = -1,\qquad
{\dot w}_2=-w_1w_2,\qquad
q\dot{w}_1+\ddot x=0.
\]
If we set $w_1(0)={\dot\theta}_0$ and $w_2(0) =
\exp(-\theta_0){\dot\phi}_0$, the solution of the above equations is
\begin{eqnarray*}
&& x(t) = -\onehalf\frac{qt^2}{q^2-1}+{\dot x}_0t+x_0,\qquad
w_1(t) =
\frac{t}{q^2-1} +{\dot\theta}_0,\\
&& w_2(t) = {\dot\phi}_0\exp\left(-\theta_0-
\frac{1}{q^2-1}(\onehalf t^2-{\dot\theta}_0t+{\dot\theta}_0q^2t)\right).
\end{eqnarray*}

Clearly, $x(t)$ has the desired form. From the reduced solution
$\check v(t)=(x(t),\dot x(t),w_1(t),w_2(t))$, we could determine the
remaining coordinates $(\theta(t),\phi(t))$ directly from the
relations $\dot\theta =w_1$ and
$\exp(-\theta)\dot\phi-\dot\theta=w_2$. The reconstruction process
as described above splits this calculation into two steps:\ first we
 calculate the horizontal lift of $\check v(t)$, and then we use it in the
reconstruction equation. In this way, we will see the effect of
changing the connection in the reconstruction equation.

The coordinates of the horizontal lift using the generalized
mechanical connection are $\hlift{\check v}(t)=(x(t),
\hlift\phi(t),\hlift\theta(t),\dot x(t),w_1(t),w_2(t))$, with respect
to the invariant basis.  They can be determined by the relation
\[
0= \mech\Omega( \dot{\hlift{\check v}}(t)) =
\mech\omega(\hlift{\check v}, T\tau\circ{\dot{\hlift{\check v}}}) =
\mech\omega\left(\hlift{\check v}, \frac{d}{dt}(\tau\circ\hlift{\check
v})\right).
\]
This equation for $\tau\circ\hlift{\check v}(t) = (x(t),
\hlift\phi(t),\hlift\theta(t))$ is
\[
\hlift{\dot\theta} = - q\dot x , \qquad \hlift{\dot\phi} = -q \dot x
\hlift\phi +q \dot x \hlift\phi=0.
\]
Therefore $\hlift\theta (t)= -q x(t)+qx_0 + \theta_0 $ and
$\hlift\phi(t)=\phi_0$.

Now we determine the curve $g(t) = (\theta_1(t),\phi_1(t))$ in $G$
such that $v=g \hlift{\check v}$ is the solution of the
Euler-Lagrange equations with the given initial values. This curve
is the solution through the identity of the reconstruction
equation $g^{-1}{\dot g} = \mech\Omega(\Gamma\circ \hlift{\check v})$.
The left hand side is
${\dot\theta}_1 E_1 + \exp(-\theta_1){\dot\phi}_1E_2$. The right hand side is
$\mech\omega(\hlift{\check v})$, or
$(w_1 + q{\dot x}) \clift{\tilde E}_1\circ \hlift{\check v}
+ (w_2 + q\phi_m {\dot x}) \clift{\tilde E}_1\circ \hlift{\check v}$.
The reconstruction equations are therefore
 \[
{\dot\theta}_1=w_1 + q {\dot x}, \qquad \exp(-\theta_1){\dot\phi}_1
= - \hlift\phi w_1 + \exp(\theta_m) w_2 - q\hlift\phi \dot x.
 \]
Solving the above equations for $(\theta_1,\phi_1)$ gives
\[
\theta_1(t) = -\onehalf t^2+(q{\dot x}_0+{\dot\theta}_0)t, \qquad
\phi_1(t) = {\dot\phi}_0t +\phi_0 \Big(1-\exp(\onehalf(2q{\dot
x}_0-t+2{\dot\theta}_0)t)\Big) .
\]
The final solution is therefore indeed
\[
\theta(t) = \theta_1(t) + \theta_m(t) =
\onehalf\frac{t^2}{q^2-1}+{\dot\theta}_0t+\theta_0,
\quad \phi(t) =
\exp(\theta_1(t))\phi_m(t)+\phi_1(t)= {\dot\phi}_0t+\phi_0.
\]

If we use the vertical lift $\Omega$ of the trivial (principal)
connection $\omega$, the equation that determines the horizontal
lift $\hlift{\check v}$ (again with group coordinates
$(\hlift\phi,\hlift\theta)$) is
$\omega(\frac{d}{dt}(\tau\circ\hlift{\check v}))=0$, or
\[
\hlift{\dot\theta} = 0, \qquad \hlift{\dot\phi} = 0,
\]
from which $\hlift\theta(t)= \theta_0$ and $\hlift\phi(t)=\phi_0$.
The equation $g^{-1}{\dot g} = \Omega(\Gamma\circ \hlift{\check v})
= \omega(\hlift{\check v})$ for $g=(\phi_1(t),\theta_1(t))$ is then
 \[
{\dot\theta}_1=w_1, \qquad \exp(-\theta_1){\dot\phi}_1 = -
\hlift\phi w_1 + \exp(\hlift\theta) w_2.
 \]
Its solution is
\[
\theta_1(t) = \onehalf \frac{t^2}{q^2-1}+{\dot\theta}_0t, \qquad
\phi_1(t) = \phi_0(1-\exp(\onehalf
\frac{t^2-2{\dot\theta}_0t+2{\dot\theta}_0q^2t}{q^2-1}))+{\dot\phi}_0t,
\]
which leads again to the same solution $(\theta(t),\phi(t))$.

The solution using the vertical lift of $\omega$ is somewhat simpler,
but this is only to be expected since $\omega$ is trivial.

\section{Conclusions and Outlook}

We have considered regular Lagrangians that are invariant under a
symmetry Lie group and we have derived the reduced Euler-Lagrange
equations, the so-called Lagrange-Poincar\'e equations. Our
framework relied on the associated Euler-Lagrange vector field and
its quotient field, rather than on the variational formalism and on
the use of well-chosen quasi-velocities. Given an integral curve of
the reduced vector field, we have shown how to reconstruct an
integral curve of the original Euler-Lagrange field by means of a
principal connection that is natural associated to the Lagrangian,
the so-called generalized mechanical connection.

In forthcoming papers, we will apply the same technique also to the
context of a different but related Lagrangian reduction technique
\cite{Routhpaper} and to the characterization of relative equilibria
\cite{releqpaper}.

\subsubsection*{Acknowledgements}

The first author is currently a Research Fellow at The University of
Michigan through a Marie Curie Fellowship.  He is grateful to the
Department of Mathematics for its hospitality.  He also acknowledges a
research grant (Krediet aan Navorsers) from the Fund for Scientific
Research - Flanders (FWO-Vlaanderen), where he is an Honorary
Postdoctoral Fellow.

The second author is a Guest Professor at Ghent University:\ he is
grateful to the Department of Mathematical Physics and Astronomy at
Ghent for its hospitality.

\end{document}